\documentclass
{amsproc}
\usepackage{amscd}

\newcommand{\coh}{cohomology}
\newcommand{\g}{Gerstenhaber}
\newcommand{\h}{Hoch\-schild}
\newcommand{\hkr}{\h-Kostant-Rosenberg}
\newcommand{\Hom}{\operatorname{Hom}}
\newcommand{\id}{\operatorname{id}}
\newcommand{\im}{\operatorname{Im}}
\newcommand{\Ker}{\operatorname{Ker}}
\newcommand{\qis}{quasi-iso\-mor\-phic}
\newcommand{\sgn}{\operatorname{sgn}}
\newcommand{\Sh}{\operatorname{Sh}}
\newcommand{\nr}{\mathbb{R}}
\newcommand{\nz}{\mathbb{Z}}

\newtheorem{theorem}{Theorem}[section]
\newtheorem{conjecture}[theorem]{Conjecture}

\theoremstyle{remark}
\newtheorem{remark}[theorem]{Remark}
\newtheorem*{ack}{Acknowledgment}

\numberwithin{equation}{section}

\begin{document}

\title{Quantizing Poisson Manifolds}
\author{Alexander A. Voronov}
\address{Department of Mathematics, M.I.T.,
2-270, 77 Massachusetts Ave., Cambridge, MA 02139-4307}
%\curraddr{}
\email{voronov@math.mit.edu}
\thanks{Research supported in part by an AMS Centennial Fellowship and 
NSF grant \#DMS-9402076.}

\date{January 15, 1997}

\begin{abstract}
This paper extends Kontsevich's ideas on quantizing Poisson
manifolds. A new differential is added to the Hodge decomposition of
the \h\ complex, so that it becomes a bicomplex, even more similar to
the classical Hodge theory for complex manifolds.
\end{abstract}

\maketitle

These notes grew out of the author's attempt to understand
Kontsevich's ideas \cite{kon:ihes} on quantizing Poisson manifolds. We
introduce a new differential on the \h\ complex, so that it becomes a
bicomplex, see Theorem~\ref{main}. This differential respects the
Hodge decomposition of the \h\ complex of a commutative algebra
discovered by \g-Schack \cite{gs:1987}. Thus, the \h\ complex becomes
similar to the $\partial$-$\bar \partial$-complex in complex
geometry. Hopefully, Hodge-theoretic ideas \emph{\`a la}
Deligne-Griffiths-Morgan-Sullivan \cite{dgms,s} will eventually result
in proving Kontsevich's Formality Conjecture, which implies local
quantization of an arbitrary Poisson manifold, a hard problem that has
been around for almost twenty years \cite{bffls}, see \cite{weinstein}
for the most state-of-the-art survey of this subject.

\begin{ack}
I am very grateful to Maxim Kontsevich, from whom I learned at least
two thirds of what is discussed in this paper. I also thank
J.-L. Brylinski, M.~Flato, M.~Gerstenhaber, J.-L. Loday, J.~Stasheff,
J.~Millson, D.~Sullivan, B.~Tsygan, and A.~Weinstein for helpful
conversations. I express my gratitude to Lew Coburn for his invitation
to participate in the 1996 Joint Summer Research Conference on
Quantization at Mount Holyoke, where this paper was delivered.
\end{ack}

\section{Kontsevich's Formality Conjecture}

\subsection{Some formalities}
\label{bracket}

Let $A = C^\infty (X)$ be the algebra of smooth functions on a smooth
real manifold $X$. Let $C^\bullet (A,A)$ be the (\emph{local\textup{)}
\h\ complex} of the algebra $A$ over $X$, \emph{i.e}., $C^n (A,A) = \{
\phi \in \Hom(A^{\otimes n}, A) \; | \; \phi(f_1,
\dots, f_n)$ is a differential operator in each entry $f_1, \dots,
f_n\}$. The \hkr\ Theorem \cite{hkr} provides the
computation of the corresponding \h\ \coh\:
\[
H^\bullet (A, A) = \Lambda^\bullet TX,
\]
which is nothing but the smooth multivector fields on $X$.

Notice that both the \h\ complex and its \coh\ (more precisely, the
suspensions thereof) are differential graded Lie algebras
(DGLA's). The complex $C^\bullet = C^\bullet(A,A)[1]$, where
$K^\bullet [1]$ denotes the suspension $K^n [1] = K^{n+1}$, $n \in
\nz$, of a complex $K^\bullet$, carries a
\emph{\g\ bracket} \cite{gerst}, which may be defined naturally, see
\cite{jim:bracket}, by observing that the \h\ cochains are exactly the
coderivations of the tensor coalgebra $T(A) =
\bigoplus_{n\ge 0} A^{\otimes n}$; then the \g\ bracket is just the
bracket of coderivations. This bracket defines a DGLA structure on
$C^\bullet$. The (suspended) $\nz$-graded vector space $H^\bullet =
\Lambda^\bullet TX[1]$ of multivector fields is a DGLA with
respect to the trivial differential $d=0$ and the
\emph{Schouten-Nijenhuis bracket} of multivector fields:
\begin{multline*}
[v_1 \wedge \dots \wedge v_m , w_1 \wedge \dots \wedge w_n]\\ =
\sum_{\substack{1 \le i \le m\\ 1\le j \le n}} (-1)^{m+i+j-1}
[v_i,w_j] \wedge v_1 \wedge \dots \wedge \widehat v_i \wedge \dots
\wedge v_m \wedge w_1 \wedge \dots \wedge \widehat w_j \wedge \dots
\wedge w_n.
\end{multline*}
Every DGLA $L^\bullet$ induces an obvious DGLA structure on its \coh\
$H^\bullet(L^\bullet)$ with the trivial differential. In this sense
the \hkr\ Theorem may be refined by saying that the \emph{\coh\ DGLA
of the \h\ complex is isomorphic to the DGLA of multivector fields},
see \cite{gersts}. Kontsevich's Formality Conjecture \cite{kon:ihes}
suggests a further, profound refinement of the \hkr\ Theorem.

\begin{conjecture}[Kontsevich's Formality Conjecture]
The \h\ complex $C^\bullet$ is \qis\ as a DGLA to its \coh\
$H^\bullet$.
\end{conjecture}
We recall that two DGLA's $L$ and $L'$ are \emph{\qis}, if there is a
chain $L=L_1 \to L_2 \leftarrow L_3 \to \dots \leftarrow L_n = L'$ of
DGLA homomorphisms all of which induce isomorphisms of \coh. Perhaps,
in this conjecture one should consider a weaker notion of
quasi-isomorphism, where the intermediate steps $L_2, L_3, \dots,
L_{n-1}$ are $L_\infty$-algebras rather than DG Lie.

\begin{remark}
There exists a natural embedding $H^\bullet \to C^\bullet$, ``a
multivector field is considered as a multiderivation of the algebra
$A$ of functions'', which induces an isomorphism of \coh. This
embedding does not satisfy the conditions of the conjecture, because
it does not respect the brackets. It is not hard to come up with a
counterexample. In rational homotopy theory, there is a similar
discouragement: the mapping $H^\bullet(X) \to \Omega^\bullet(X)$ which
takes a \coh\ class to its harmonic representative is a
quasi-isomorphism, but the product of two harmonic forms is not
harmonic in general. Nevertheless, the two differential graded
associative algebras $H^\bullet(X)$ and $\Omega^\bullet(X)$ are \qis\
for a compact K\"ahler $X$, see Section~\ref{rht}.
\end{remark}

We will take the ``physical'' point of view and discuss evidence for
the Formality Conjecture after seeing what implications it has.

\subsection{Deformation quantization of Poisson manifolds}

\begin{sloppypar}
Recall that a \emph{deformation quantization} \cite{bffls} of a
Poisson manifold $X$, whose algebra of smooth functions will be
denoted by $A$, as above, is a formal deformation of $A$ in the
direction of the Poisson bracket. More precisely, it is a
multiplication $a*b$ on $A[[h]] = A \otimes \nr[[h]]$ making it an
associative $\nr[[h]]$-algebra, such that for $a,b \in A$
\[
a * b = ab + \{a,b\}h + B(a,b) h^2 + \dots,
\]
where $ab$ is the usual, undeformed multiplication and $\{a,b\}$ is
the Poisson bracket. When the Poisson bracket is nondegenerate,
\emph{i.e}., coming from a symplectic structure, the existence of
deformation quantization was proven by De Wilde and Le\-comte
\cite{dwlc} and Fedosov \cite{fedosov}. When the Poisson bracket is
arbitrary, the existence of deformation quantization (even locally,
for $\nr^n$) is an open problem.
\end{sloppypar}

The remarkable fact noticed by Kontsevich is that if you assume the
Formality Conjecture, the problem of quantization will be solved.
\begin{theorem}[Kontsevich]
The Formality Conjecture for a manifold $X$ implies deformation
quantization of any Poisson structure on $X$.
\end{theorem}

\begin{proof}
We will only sketch the idea of the proof; (some) details may be found
in Kontsevich's Berkeley lectures \cite{kon:def}.

According to Deligne-Schlessinger-Stasheff-Goldman-Millson's approach
to deformation theory, see \cite{gm,millson,ss}, with each DGLA
$L^\bullet$ one can associate the \emph{formal moduli space} $M= \{
\gamma \in L^1 \; | \; d\gamma + \frac{1}{2} [\gamma, \gamma] =
0\}/\exp (L^0)$, where $\exp (L^0)$ is the Lie group corresponding to
the Lie algebra $L^0$, and the defining equation $d \gamma +
\frac{1}{2} [\gamma, \gamma] = 0$ is called the \emph{Maurer-Cartan
equation}. If $L_1^\bullet \to L_2^\bullet$ is a quasi-isomorphism of
DGLA's, then the corresponding formal moduli spaces can be
identified. This is done using the standard machinery of minimal
models.

Formal deformations are usually formal paths in the formal moduli
spaces. Consider the cases of the above two DGLA's associated to a
manifold $X$. The formal moduli space associated to the DGLA
$C^\bullet$ is $M_Q = \{\gamma \in C^2(A,A) \; | \; d\gamma +
\frac{1}{2} [\gamma,\gamma] = 0\} /\operatorname{GL}(A)$. Since the
\h\ differential $d$ is equal to the bracket $[m,\,]$ with the
multiplication two-cocycle $m \in C^2(A,A)$, the Maurer-Cartan
equation may be rewritten as $[m+\gamma, m+ \gamma] = 0$, which is
equivalent to the associativity of $m + \gamma \in \Hom(A \otimes A,
A)$ understood as a new multiplication. A deformation quantization (in
an arbitrary direction) is then a formal path originating at $\gamma =
0$ in the formal moduli space $M_Q$. If this deformation quantizes a
Poisson bracket, the tangent vector to the formal path should coincide
with the Poisson bracket.

The formal moduli space associated to the DGLA $H^\bullet$ of
multivector fields on $X$ is $M_P = \{\gamma \in \Lambda^2 TX \; | \;
[\gamma, \gamma] = 0\} / \exp (\operatorname{Vect} X)$. The
Maurer-Cartan equation in this case is equivalent to the Jacobi
identity for the skew bracket $\{f,g\} = \gamma (df,dg)$ of functions
on $X$. Thus a solution of the Maurer-Cartan equation is a Poisson
structure on $X$.

Now suppose that the Formality Conjecture is true. Then the two moduli
spaces $M_Q$ and $M_P$ are identified. Given a Poisson structure on
$X$, we can connect it by a straight line with the origin in the
moduli space of Poisson structures. Consider this line as a formal
path. Using the isomorphism of the moduli spaces, we have a formal
path in the moduli space of quantizations, which is a deformation
quantization we were looking for.
\end{proof}

\subsection{Evidence for the Formality Conjecture}

It is known that every nondegenerate Poisson structure can be
quantized \cite{dwlc,fedosov}. Moreover, the following analogue of the
conjecture related to the nondegenerate case is true. Consider the \h\
DGLA of the function algebra with respect to the deformed
multiplication on a symplectic manifold $X$. The other DGLA will be
the multivector fields on $X$ with the differential being the
Schouten-Nijenhuis bracket with the canonical Poisson tensor on
$X$. Then the two DGLA's are \qis. One uses Fedosov's connection to
prove this fact, \cite{kon:ihes}.

Different evidence comes from quantizing an arbitrary Poisson
structure in the Lie-theoretic context. The recent theorem of
P.~Etingof and D.~Kazhdan \cite{ek} solves the conjecture of Drinfeld
asserting that every Poisson Lie group has a canonical quantization.

Mirror Symmetry predicts that for a Calabi-Yau manifold $Y$, the
corresponding holomorphic version $R\Gamma(Y,C^\bullet)$ of the \h\
DGLA for the sheaf of holomorphic functions on $Y$ gives rise to a
smooth formal moduli space, which may be interpreted as the moduli
space of ``noncommutative Calabi-Yau manifolds''. On the other hand,
one can show that the holomorphic multivector field DGLA $R\Gamma
(Y,H^\bullet)$ produces a smooth formal moduli space. Thus, if
$C^\bullet$ and $H^\bullet$ were known to be \qis, it would prove the
smoothness of the first moduli space as confirmed by Mirror Symmetry.

\begin{sloppypar}
\subsection{Formality in rational homotopy theory}
\label{rht}

Kontsevich's Formality Conjecture has a very close analogy with the
Deligne-Griffiths-Morgan-Sullivan Formality Theorem \cite{dgms}:
\emph{the Sullivan model of a compact K\"ahler manifold $X$ is
formal}. The \emph{Sullivan model} of $X$ may be represented by the
differential graded commutative algebra (DGA) $\Omega^\bullet (X)$ of
smooth differential forms on $X$. Formality means that $\Omega^\bullet
(X)$ is \qis\ to its \coh\ DGA $H^\bullet(X)$. A simple way to prove
this is using Hodge theory, see \cite{s}: decompose the de Rham
differential into the holomorphic and antiholomorphic parts: $d =
\partial + \bar \partial$. Standard Hodge-theoretic arguments (the
$\partial$-$\bar \partial$-Lemma of \cite{dgms}) imply that $(\Ker
\partial, d) \subset (\Omega^\bullet (X), d)$ is an embedding of
DGA's, which is a quasi-isomorphism. On the other hand, the natural
morphism $(\Ker \partial, d) \to (\Ker \partial / \im \partial, d) =
(H^\bullet(X), 0)$ of DGA's is also a quasi-isomorphism for the same
reasons.
\end{sloppypar}

\section{Hodge theory for the \h\ complex}

In this section, we are going to develop Hodge theory in the \h\
context. The construction of Hodge decomposition of the \h\ complex of
a commutative algebra $A$ over a field of characteristic zero goes
back to \g\ and Schack \cite{gs:1987}, who decomposed the
Hochschild complex $C^\bullet(A,A)$ into the direct sum of
$C^{p,q}(A,A)$, with the \h\ differential $d$ acting like $\bar
\partial$ in the Dolbeault complex: $d: C^{p,q} (A,A) \to C^{p,q+1}
(A,A)$. Here we add a new ingredient to \g-Schack's Hodge theory: we
define an extra, $\partial$-like differential $d': C^{p,q} (A,A) \to
C^{p-1,q} (A,A)$ on the \h\ complex, so that it becomes a
bicomplex. This bicomplex is similar to the $\partial$-$\bar
\partial$-complex of a compact K\"ahler manifold: the total
cohomology of the bicomplex is equal to the cohomology of one of the
differentials. Our new differential is also similar to the
differential $B$ of the cyclic \coh\ complex. Together with the \h\
differential, the differential $B$ provides the cyclic cohomology
complex with the structure of a bicomplex and, moreover, respects the
Hodge decomposition of the cyclic \coh\ complex in a similar way, see
J.-L. Loday \cite{loday:hodge}. Another similarity between the cyclic
$B$ and our differential is that the \coh\ of both vanish.

We will recall Hodge decomposition of the \h\ complex, following the
modification of M.~Ronco, A.~B. Sletsj{\o}e, and H.~L. Wolfgang, see
\cite{bw} for more detail. Let $r$ and $s$ be positive integers and $n
= r+s$. The \emph{shuffle product} of tensors $a_1 \otimes \dots
\otimes a_{r} \in A^{\otimes r}$ and $a_{r+1} \otimes \dots \otimes
a_{n} \in A^{\otimes s}$ is the element
\[
\sum \sgn (\sigma) a_{\sigma^{-1} (1)} \otimes \dots \otimes 
a_{\sigma^{-1} (n)} \in A^{\otimes n},
\]
where the summation runs over those $\sigma \in S_n$ for which
$\sigma(1) < \sigma (2) < \dots < \sigma (r)$ and $\sigma(r+1) <
\dots < \sigma (n)$. Let $\Sh^k$ denote the image of shuffle products
of $k$ elements in the tensor algebra $T(A) = \bigoplus_{n \ge 0}
A^{\otimes n}$. By definition $\Sh^0 = T(A)$ and $\Sh^1 = \bigoplus_{n
> 0} A^{\otimes n}$. We have a filtration of the tensor algebra
\[
T(A) = \Sh^0 \supset \Sh^1 \supset \Sh^2 \supset \Sh^3 \dots
\]
Define
\[
C^{p,q} = \Hom(\Sh^p \cap A^{\otimes p+q} / \Sh^{p+1} \cap A^{\otimes
p+q}, A),
\]
where $p,q \ge 0$. Of course, one can describe $C^{p,q}$ as $A$-valued
functionals $\phi$ on the subspace of $A^{\otimes p+q}$ generated by
the shuffle products $\Sh^p$ of $p$ elements, such that $\phi$
vanishes on the shuffle products of $p+1$ elements. One can check that
the \h\ differential induces a mapping $d:C^{p,q} \to C^{p,q+1}$, $d^2
=0$, and $C^n(A,A) \cong \bigoplus_{p+q = n} C^{p,q}$.

This gives the \emph{Hodge decomposition} of the \h\ complex into the
direct sum of complexes:
\begin{gather*}
\dots\\
\begin{CD}
0 @>>> C^{2,0} @>d>> C^{2,1} @>d>> C^{2,2} @>d>> \dots\\
0 @>>> C^{1,0} @>d>> C^{1,1} @>d>> C^{1,2} @>d>> \dots\\
0 @>>> C^{0,0} @>d>> 0 @>>> 0 @>>> \dots
\end{CD}
\end{gather*}
The second-to-last row is known as the \emph{Harrison
complex}\footnote{The last row does not have a name yet, but will
hopefully acquire one soon.}. The Hodge decomposition of the \h\
complex induces one on the \h\ \coh: $H^n(A,A) \cong \bigoplus_{p+q =
n} H^{p,q}(A,A)$.

\begin{theorem}
\label{main}
\begin{enumerate}
\item There exists a differential $d':C^{p,q} \to C^{p-1,q}$ which
is a derivation of the \g\ bracket and defines the structure of a
bicomplex on the \h\ complex:
\begin{gather*}
\dots\\
\begin{CD}
@VVV @VVd'V @VVd'V @VVd'V\\
0 @>>> C^{2,0} @>d>> C^{2,1} @>d>> C^{2,2} @>d>> \dots\\
@VVV @VVd'V @VVd'V @VVd'V\\
0 @>>> C^{1,0} @>d>> C^{1,1} @>d>> C^{1,2} @>d>> \dots\\
@VVV @VVd'V @VVd'V @VVd'V\\
0 @>>> C^{0,0} @>d>> 0 @>>> 0 @>>> \dots\\
@VVV @VVd'V @VVV @VVV\\
0 @>>> 0  @>>> 0 @>>> 0 @>>> \dots
\end{CD}
\end{gather*}

\item The spectral sequence associated to the first filtration
$'F^p = \bigoplus_{i\le p} C^{i,j}$ is convergent: $'E_\infty =
H^\bullet (C^{\bullet, \bullet}, d+d')$. Moreover, $'E_1$ is equal to
the \h\ \coh\ $H^\bullet (A,A)$.

\item The \coh\ of the differential $d'$ vanishes. The spectral
sequence associated to the second filtration $''F_q = \bigoplus_{j\ge
q} C^{i,j}$ collapses at $''E_1$, which is equal to 0.

\item Suppose that $A$ is the algebra of smooth functions on a manifold or
regular functions on a nonsingular affine scheme. Then the first
spectral sequence collapses at $'E_1$, which is equal to $H^{\bullet,
0} (A,A)$, the space of global multivector fields. This coincides with
the total \coh\ of the bicomplex.
\end{enumerate}
\end{theorem}

\begin{remark}
The differential $d'$, being a derivation of the \g\ brack\-et of
degree $-1$, defines the structure of a DGLA on the \h\ complex
$(C^\bullet(A,A)[1], d')$. However, the total complex
$(C^\bullet(A,A)[1], d+d')$ is only a differential $\nz/2\nz$-graded
Lie algebra: the degree of the total differential $d+d'$ is equal to
one modulo two.
\end{remark}

\begin{proof}
1. Define $d'$ as the inner derivation
\[
d' \phi = [1,\phi],
\]
where $[\,,]$ is the \g\ bracket (see Section~\ref{bracket}) and $1
\in A = C^0(A,A)$ is the unit element of $A$. 
In other words,
\begin{multline*}
(d'\phi) (a_1, \dots, a_{n-1}) \\
= \phi(a_1, \dots, a_{n-1}, 1) - \phi(a_1, \dots, a_{n-2}, 1,
a_{n-1}) + \dots +(-1)^{n-1} \phi(1, a_1, \dots, a_{n-1}).
\end{multline*}
Then $(d')^2 =0$, because $[1,1]=0$, and $d'd + dd'=0$, because $d =
[m,\,]$, where $m \in C^2(A,A)$ is the multiplication cocycle, and
$[m,1] = 0$: $[m,1] (a) = m(1,a) - m(a,1) = 1\cdot a - a \cdot 1 = 0$
for any $a \in A$. It is also clear that $d': C^n(A,A) \to C^{n-1}
(A,A)$. Let us verify that moreover $d': C^{p,q} \to
C^{p-1,q}$. Indeed, if $\phi$ is an $A$-valued functional defined on
the $p$-shuffles $\Sh^p$, then $d'\phi$ is obviously defined on
$p-1$-shuffles via the natural mapping $\Sh^{p-1} \to \Sh^p$, the
shuffle product with $1 \in A$. For the same reason, if $\phi$
vanishes on $\Sh^{p+1}$, then $d'\phi$ will vanish on $\Sh^p$. This
proves the first statement of the theorem.

2. The first spectral sequence is convergent, because the first
filtration is regular: moreover, $'F^{-1}=0$. Since $'F^p/'F^{p+1} =
C^{p, \bullet}$ with the differential $d$, $'E^p_1 = H^\bullet
('F^p/'F^{p+1}, d) = H^{p,\bullet}(A,A)$.

3. We will define for each $n \ge 0$ a null-homotopy $k: C^n(A,A) \to
C^{n+1}(A,A)$ on the complex $C^{\bullet}(A,A)$ with respect to the
differential $d'$:
\[
(k\phi)(a_1, \dots, a_{n+1}) = \frac{1}{n+1}\sum_{i=1}^{n+1}
(-1)^{i-1} a_i \phi(a_1, \dots, \widehat a_i, \dots, a_{n+1}).
\]
If $\phi = 0$ on $\Sh^{r}$, then $k\phi = 0$ on $\Sh^{r+1}$, therefore
$k$ is well-defined on $C^{p,q}$ and maps it to $C^{p+1,q}$. A
straightforward computation shows that $kd' + d'k = \id$. Thus, the
\coh\ of $d'$ vanishes. Since $''F_q/''F_{q+1} = C^{\bullet, q}$ with the
differential $d'$, $''E^q_1 = H^\bullet (''F_q/''F_{q+1}, d') = 0$,
and the second spectral sequence collapses.

4. If $A$ is a regular algebra of functions, its \h\ \coh\ $H^\bullet
(A,A)$ is equal to the space of multivector fields, see
\cite{hkr}. The multivector fields are skew multiderivations of $A$ and
therefore project bijectively on $H^{\bullet,0}(A,A)$. In this case,
the differential $d'$ vanishes on all \h\ cocycles, because
derivations of $A$ vanish on constants. Therefore $'E_1 = 'E_2 = \dots
= 'E_\infty = H^{\bullet,0}(A,A)$. The computation of the total \coh\
then follows from Part 2.
\end{proof}

\bibliographystyle
{amsalpha}
%{/usr/local/lib/texmf/tex/latex/packages/amslatex/classes/amsalpha}
%\bibliography{../appl/appl,../gregg/op}

\begin{thebibliography}{DGMS75}

\bibitem[BFF{\etalchar{+}}78]{bffls}
F.~Bayen, M.~Flato, C.~Fronsdal, A.~Lichnerowicz, and D.~Sternheimer,
  \emph{Deformation theory and quantization. {I}. {D}eformations of symplectic
  structures}, Ann. Physics \textbf{111} (1978), no.~1, 61--110.

\bibitem[BW95]{bw}
N.~Bergeron and H.~L. Wolfgang, \emph{The decomposition of {H}ochschild
  cohomology and {G}erstenhaber operations}, J. Pure Appl. Algebra \textbf{104}
  (1995), no.~3, 243--265.

\bibitem[DGMS75]{dgms}
P.~Deligne, P.~Griffiths, J.~Morgan, and D.~Sullivan, \emph{Real homotopy
  theory of {K}\"ahler manifolds}, Invent. Math. \textbf{29} (1975), no.~3,
  245--274.

\bibitem[DWL83]{dwlc}
M.~De~Wilde and P.~B.~A. Lecomte, \emph{Existence of star-products and of
  formal deformations of the {P}oisson {L}ie algebra of arbitrary symplectic
  manifolds}, Lett. Math. Phys. \textbf{7} (1983), no.~6, 487--496.

\bibitem[EK96]{ek}
P.~Etingof and D.~Kazhdan, \emph{Quantization of {L}ie bialgebras. {I}},
  Selecta Math. (N.S.) \textbf{2} (1996), no.~1, 1--41.

\bibitem[Fed85]{fedosov}
B.~V. Fedosov, \emph{Formal quantization}, Some problems in modern mathematics
  and their applications to problems in mathematical physics (Russian)
  (Moscow), Moskov. Fiz.-Tekhn. Inst., Moscow, 1985, pp.~129--136, {\rm vi}.

\bibitem[Ger63]{gerst}
M.~Gerstenhaber, \emph{The cohomology structure of an associative ring}, Ann.
  of Math. \textbf{78} (1963), 267--288.

\bibitem[GM90]{gm}
W.~M. Goldman and J.~J. Millson, \emph{The homotopy invariance of the
  {K}uranishi space}, Illinois J. Math. \textbf{34} (1990), no.~2, 337--367.

\bibitem[GS87]{gs:1987}
M.~Gerstenhaber and S.~D. Schack, \emph{A {H}odge-type decomposition for
  commutative algebra cohomology}, J. Pure and Appl. Alg. \textbf{48} (1987),
  229--247.

\bibitem[GS88]{gersts}
M.~Gerstenhaber and S.~D. Schack, \emph{Algebraic cohomology and deformation
  theory}, Deformation Theory of Algebras and Structures and Applications
  (M.~Hazewinkel and M.~Gerstenhaber, eds.), Kluwer Academic Publishers, 1988,
  pp.~11--264.

\bibitem[HKR62]{hkr}
G.~Hochschild, B.~Kostant, and A.~Rosenberg, \emph{Differential forms on
  regular affine algebras}, Trans. Amer. Math. Soc. \textbf{102} (1962),
  383--408.

\bibitem[Kon95a]{kon:ihes}
M.~Kontsevich, \emph{Lectures at {IHES}}, Fall 1995.

\bibitem[Kon95b]{kon:def}
M.~Kontsevich, \emph{Lectures on deformation theory}, Preprint, University of
  California, Berkeley, Spring 1995.

\bibitem[Lod89]{loday:hodge}
J.-L. Loday, \emph{Op\'erations sur l'homologie cyclique des alg\`ebres
  commutatives}, Invent. Math. \textbf{96} (1989), no.~1, 205--230.

\bibitem[Mil91]{millson}
J.~J. Millson, \emph{Rational homotopy theory and deformation problems from
  algebraic geometry}, Proceedings of the International Congress of
  Mathematicians, Vol.\ I, II (Kyoto, 1990) (Tokyo), Math. Soc. Japan, 1991,
  pp.~549--558.

\bibitem[SS85]{ss}
M.~Schlessinger and J.~Stasheff, \emph{The {L}ie algebra structure of tangent
  cohomology and deformation theory}, J. Pure Appl. Algebra \textbf{38} (1985),
  no.~2-3, 313--322.

\bibitem[Sta93]{jim:bracket}
J.~Stasheff, \emph{The intrinsic bracket on the deformation complex of an
  associative algebra}, J. Pure Appl. Algebra \textbf{89} (1993), 231--235.

\bibitem[Sul77]{s}
D.~Sullivan, \emph{Infinitesimal computations in topology}, Inst. Hautes
  \'Etudes Sci. Publ. Math. (1977), no.~47, 269--331 (1978).

\bibitem[Wei95]{weinstein}
A.~Weinstein, \emph{Deformation quantization}, Ast\'erisque (1995), no.~227,
  Exp.\ No.\ 789, 5, 389--409, S\'eminaire Bourbaki, Vol.\ 1993/94.

\end{thebibliography}

%\end{document}

\newcommand{\etalchar}[1]{$^{#1}$}
\providecommand{\bysame}{\leavevmode\hbox to3em{\hrulefill}\thinspace}

\end{document}